\documentclass[11pt]{article}
\usepackage[T1]{fontenc}
\usepackage{theorem}
\usepackage{amsmath,amssymb,latexsym}
\usepackage[final]{graphicx}
\usepackage{ae}
\newtheorem{thm}{Theorem}[section]
\newtheorem{prop}[thm]{Proposition}
\newtheorem{lm}[thm]{Lemma}
\newtheorem{cor}[thm]{Corollary}
{\theorembodyfont{\rmfamily}

\newtheorem{remark}[thm]{Remark}}
\newcommand{\res}{\mathop{\rm res}}
\newcommand{\const}{\mathop{\rm const}}

\newcommand{\field}[1]{\mathbb{#1}}
\newcommand{\R}{\field{R}}

\newcommand{\N}{\field{N}}
\newcommand{\C}{\field{C}}

\renewcommand{\Re}{\mathop{\rm Re}}
\renewcommand{\Im}{\mathop{\rm Im}}

\newcommand{\isdef}{\stackrel{\text{\tiny def}}{=}}

\DeclareRobustCommand{\qed}{%
\ifmmode 
\else \leavevmode\unskip\penalty9999 \hbox{}\nobreak\hfill \fi
\quad\hbox{\qedsymbol}}
\newcommand{\openbox}{\leavevmode
\hbox to.77778em{%
\hfil\vrule
\vbox to.675em{\hrule width.6em\vfil\hrule}%
\vrule\hfil}}
\newcommand{\qedsymbol}{\openbox}
\newcommand{\proofname}{Proof}
\newenvironment{proof}[1][\proofname]{\par
\normalfont \trivlist \item[\hskip\labelsep   \itshape #1. ]
\ignorespaces
}{%
\qed\endtrivlist } 

\author{A. Mart\'{\i}nez-Finkelshtein,\footnote{Corresponding author. Email: \texttt{andrei@ual.es}}
\qquad J. F. S\'{a}nchez-Lara\\ University of Almer\'{\i}a, Spain}
\title{Shannon entropy of symmetric Pollaczek polynomials}%
\begin{document}
\maketitle
\begin{abstract}
We discuss the asymptotic behavior (as $n\to \infty$)  of the
entropic integrals
$$
E_n= -  \int_{-1}^1 \log \big( p^2_n(x) \big) p^2_n(x) w(x)\, d x
\,,
$$
and
$$
 F_n  =  -\int_{-1}^1
\log \left( p_n^2(x)w(x) \right ) p_n^2(x) w(x)\, dx,
$$
when $w$ is the symmetric Pollaczek weight on $[-1,1]$ with main
parameter $\lambda\geq 1$, and $p_n$ is the corresponding
orthonormal polynomial of degree $n$. It is well known that $w$ does
not belong to the Szeg\H{o} class, which implies in particular that
$E_n\to -\infty$. For this sequence we find the first two terms of
the asymptotic expansion. Furthermore, we show that $F_n \to \log
(\pi)-1$, proving that this ``universal behavior'' extends beyond
the Szeg\H{o} class. The asymptotics of $E_n$ has also a curious
interpretation in terms of the mutual energy of two relevant
sequences of measures associated with $p_n$'s.
\end{abstract}

\section{Introduction and statement of results}
Different information measures, and in particular, the Shannon
entropy, has found application in many branches of science. In
quantum mechanics, the uncertainty in the localization of a particle
in ordinary space is quantitatively measured by the so-called
position information entropy
\begin{equation*}\label{uno}
S_{\rho} = - \int \rho (\vec{r}) \log \rho (\vec{r}) d \vec{r}\,,
\end{equation*}
of the probability density $\rho (\vec{r}) = \left| \psi (\vec{r})
\right|^{2}$, where $\psi (\vec{r})$ is the wavefunction of its
dynamical state. This functional leads, for instance, to a
stronger version of the celebrated Heisenberg's uncertainty
principle, a fundamental law of nature
\cite{Bialynicki-Birula:75}. This fact and the effective
implementation of the density functional theory of complex
many-electron systems \cite{Parr89}, which uses the
single-particle density as the basic variable, are responsible for
the fact that the study of the entropy has become a standard tool
in atomic and molecular physics, and in condensed matter theories.
The exact or explicit determination of the information entropies
of complex many-particle systems is an extremely difficult
problem. Only recently a small progress has been achieved yielding
in some cases closed formulas for the information entropies of the
simplest 1-dimensional single-particle systems and the
three-dimensional systems of particles moving in a central or
spherically symmetric potential. For these systems the
wavefunctions are expressible in terms of some special functions,
and the determination of the corresponding information entropies
boils down naturally to the computation of entropic functionals
for sequences of orthogonal polynomials (cf.\
\cite{Yanez:94,Yanez:99}; a state-of-the art of this topic up to
2001 is given in \cite{Dehesa:01}). In particular, given a
positive unit itegrable weight $w$ on $[-1,1]$, and the sequence
of corresponding orthonormal polynomials $\{p_n\}_{n\geq 0}$, we
may define the Shannon entropy of these polynomials either as
\begin{equation} \label{def_E}
E_n=E_n(w ) \isdef -  \int_{-1}^1 \log \big( p^2_n(x) \big) p^2_n(x)
w(x)\, d x \,,
\end{equation}
or as
\begin{equation} \label{def_F}
 F_n =F_n(w )\isdef  -\int_{-1}^1
\log \left( p_n^2(x)w(x) \right ) p_n^2(x) w(x)\, dx.
\end{equation}
They are obviously related by $E_n(w)-F_n(w)=G_n(w)$, where
\begin{equation} \label{def_G}
G_n =G_n(w)\isdef \int_{-1}^1 \log(w(x)) p_n^2(x) w(x)\, dx.
\end{equation}
Hence, we are faced with two different problems. One is the
explicit computation of \eqref{def_E}--\eqref{def_F} for fixed
$n$'s (either as a closed formula or numerically). Observe that a
naive evaluation of these functionals by means of quadratures
encounters the difficulty of the zeros of $p_n$, that belong to
the interval of orthogonality. Some of the contributions in this
sense are
\cite{Buyarov:97,MR1790053,Buyarov04,Dehesa:97,Yanez:94}. A second
problem is the study of the asymptotic behavior of $\{E_n\}$,
$\{F_n\}$, and $\{G_n\}$ when $n \to \infty$, which has a special
interest in the analysis of the highly-excited (Rydberg) states of
numerous quantum-mechanical systems 
\cite{Yanez:94}. In this sense there have been important
contributions in the last few years \cite{Aptekarev:95,
Buyarov:99a, MR1790053, Dehesa:01, Dehesa:98, Levin03, Lara02}. In
a recent paper \cite{Beckermann04} the authors have studied the
asymptotic behavior of these functionals under the assumption that
the weight of orthogonality satisfies the Szeg\H{o} condition,
\begin{equation}\label{szego}
    \int_{-1}^1 \frac{\log\left( w(x)\right)}{\sqrt{1-x^2}}\, dx
    >-\infty\,.
\end{equation}
Under an additional assumption on the growth of the polynomials on
the interval of orthogonality they proved that both $E_n$ and $F_n$
(and in consequence, also $G_n$) converge, and
\begin{equation}\label{limitF}
    \lim_{n} F_n(w)=\log(\pi)-1
\end{equation}
(notice that $F_n$ is taken here with a slightly different
normalization than in \cite{Beckermann04}). The authors of
\cite{Beckermann04} conjectured that the limit in \eqref{limitF} is
valid for a larger class of weights; from their work it follows also
that if $w>0$ on $(-1,1)$ does not satisfy \eqref{szego}, then
$E_n(w)$ and $G_n(w)$ diverge to $-\infty$. The Pollaczek
polynomials constitute the first and the best known example of a
family of orthogonal polynomials on $[-1,1]$ with respect to a
weight \emph{not satisfying} the Szeg\H{o} condition \eqref{szego}.
In this paper we deal with the \emph{symmetric} Pollaczek
polynomials, $p_n^{\lambda}(x;a)$, that depend on two real
parameters, $ \lambda>0$, $  a\geq 0$, and that may be defined by
the recurrence relation
\begin{equation}\label{RecPoll}
x p_n^{\lambda}(x;a)=a_{n+1}\, p_{n+1}^{\lambda}(x;a ) + a_n\,
p_{n-1}^{\lambda}(x;a ), \quad p_{-1}^{\lambda}(x;a)=0\,, \quad
p_0^{\lambda}(x;a)=1\,,
\end{equation}
with the coefficients
\begin{equation}
a_n = \frac{1}{2}
\sqrt{\frac{n(n+2\lambda-1)}{(n+\lambda+a)(n+\lambda+a-1)}}\,.\label{CoefRecPoll}
\end{equation}
It is known (see \cite[Appendix]{szego:1975}) that these polynomials
are orthonormal on $ [-1,1]$ with respect to the unitary weight
function
\begin{equation} \label{WeightPollaczek}
w_{\lambda}(x;a)=\frac{2^{2\lambda} (\lambda+a)}{2
\pi~\Gamma(2\lambda)}\, (1-x^2)^{\lambda-1/2}~e^{(2\arccos
x-\pi)\frac{ax }{\sqrt{1-x^2}}}~\left|\Gamma\left(\lambda+i
 \frac{ax }{\sqrt{1-x^2}}\right)\right|^2,
\end{equation}
where $\Gamma(x)$ denotes the gamma function. From
\eqref{WeightPollaczek} it is clear that Pollaczek polynomials
$p_{n}^{\lambda}(x;0 )$ (that is, for $a=0$) reduce to orthonormal
Gegenbauer polynomials with parameter $\lambda $. In the sequel,
whenever it cannot lead us into confusion, we omit the explicit
reference to the parameters $\lambda $ and $a$ from the notation of
the polynomials.

Our main goal is to study the asymptotic behavior of the sequences
$E_n(w)$, $F_n(w)$ and $G_n(w)$ as $n\to \infty$, when $w$ is the
symmetric Pollaczek weight. We can summarize our main results saying
that limit \eqref{limitF} is proved to be valid also for
$w=w_{\lambda}(\cdot ;a)$, with the restriction $\lambda \geq 1$
(so, this fact extends beyond the Szeg\H{o} class, as conjectured),
and we find the main part of the asymptotic expansion of $E_n(w)$
and $G_n(w)$, up to the $o(1)$ terms. Namely, we establish the
following: for $F_n$ we prove that indeed, also for the symmetric
Pollaczek weight, limit \eqref{limitF} is still valid:
\begin{thm}\label{thFn}
For the symmetric Pollaczek weight $w=w_{\lambda}(\cdot ;a)$, with
$a\geq 0$ and $\lambda\geq 1$,
$$
F_n(w)=\log(\pi)-1+o(1), \qquad n \rightarrow\infty.
$$
\end{thm}
\begin{remark}
The restriction $\lambda\geq 1$ comes from the method of proof; we
believe that Theorem \ref{thFn} is valid for the whole range of
the parameter $\lambda $, that is, for $\lambda >0$.
\end{remark}
For the divergent sequence $\{G_n\}$ we find the first two terms of
its asymptotic expansion:
\begin{thm}\label{thGn}
For the symmetric Pollaczek weight $w=w_{\lambda}(\cdot ;a)$, with
$a\geq 0$ and $\lambda>0$,
\begin{align*}
G_n(w) &= -2a\log(n) +2a + \log\left(\frac{
\Gamma(\lambda+a)\Gamma(\lambda+a+1) }{\pi~\Gamma(2\lambda )}\right)
+o(1), \quad n\to \infty\,.
\end{align*}
\end{thm}
As a straightforward corollary we obtain
\begin{cor}\label{TeoAsPoll}
For the symmetric Pollaczek weight $w=w_{\lambda}(\cdot ;a)$, with
$a\geq 0$ and $\lambda\geq 1$,
\begin{equation}\label{asymptoticFormulaforE}
E_n(w)= -2a\log(n)+\tau (\lambda ;a)+o(1),\qquad n\rightarrow\infty,
\end{equation}
where
\begin{equation}\label{def_tau}
\tau (\lambda ;a)\isdef 2a-1 + \log\left(\frac{
\Gamma(\lambda+a)\Gamma(\lambda+a+1) }{ \Gamma(2\lambda )}\right)
\,.
\end{equation}
\end{cor}
\begin{remark}
The value $\tau(\lambda ;0)=-1+\log\left(\Gamma(\lambda
)\Gamma(\lambda +1)/ \Gamma(2\lambda ) \right)$ matches  $\lim_n
E_n$ for orthonormal Gegenbauer polynomials, found in
\cite{Aptekarev:95}.
\end{remark}

\begin{figure}[p]
\centering
\includegraphics[scale=0.65]{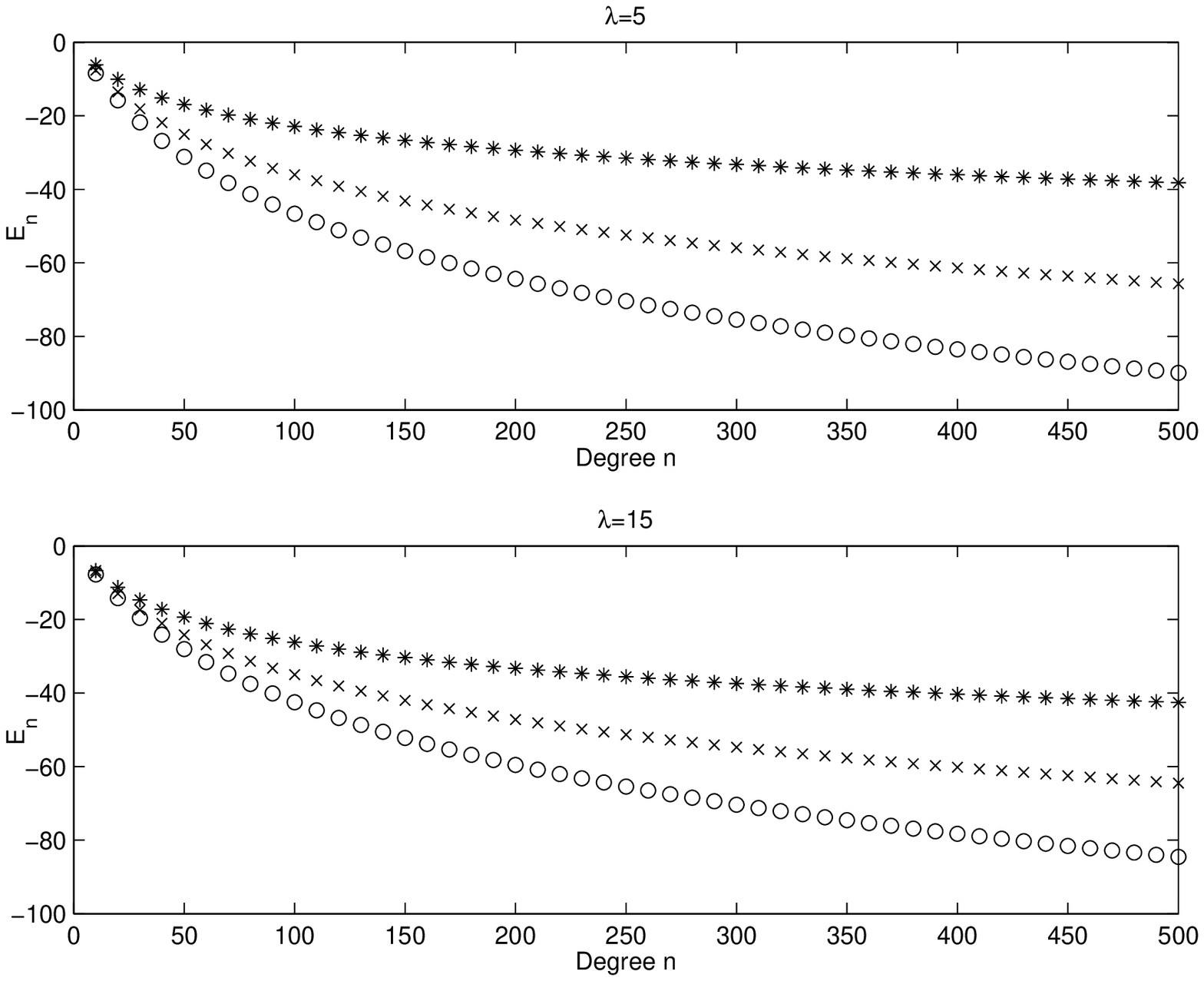}
\caption{Entropy  $E_n(w_{\lambda}(\cdot;a))$ for $n =1, 2, \dots,
500$, with $\lambda =5$ (upper) and $\lambda=15$ (lower). We use the
values $a=5$ (`\texttt{*}'), $a=10$ (`\texttt{x}') and $a=15$
(`\texttt{o}'). \label{grafPoll1}}
\end{figure}

\begin{figure}[p]
\centering
\includegraphics[scale=0.65]{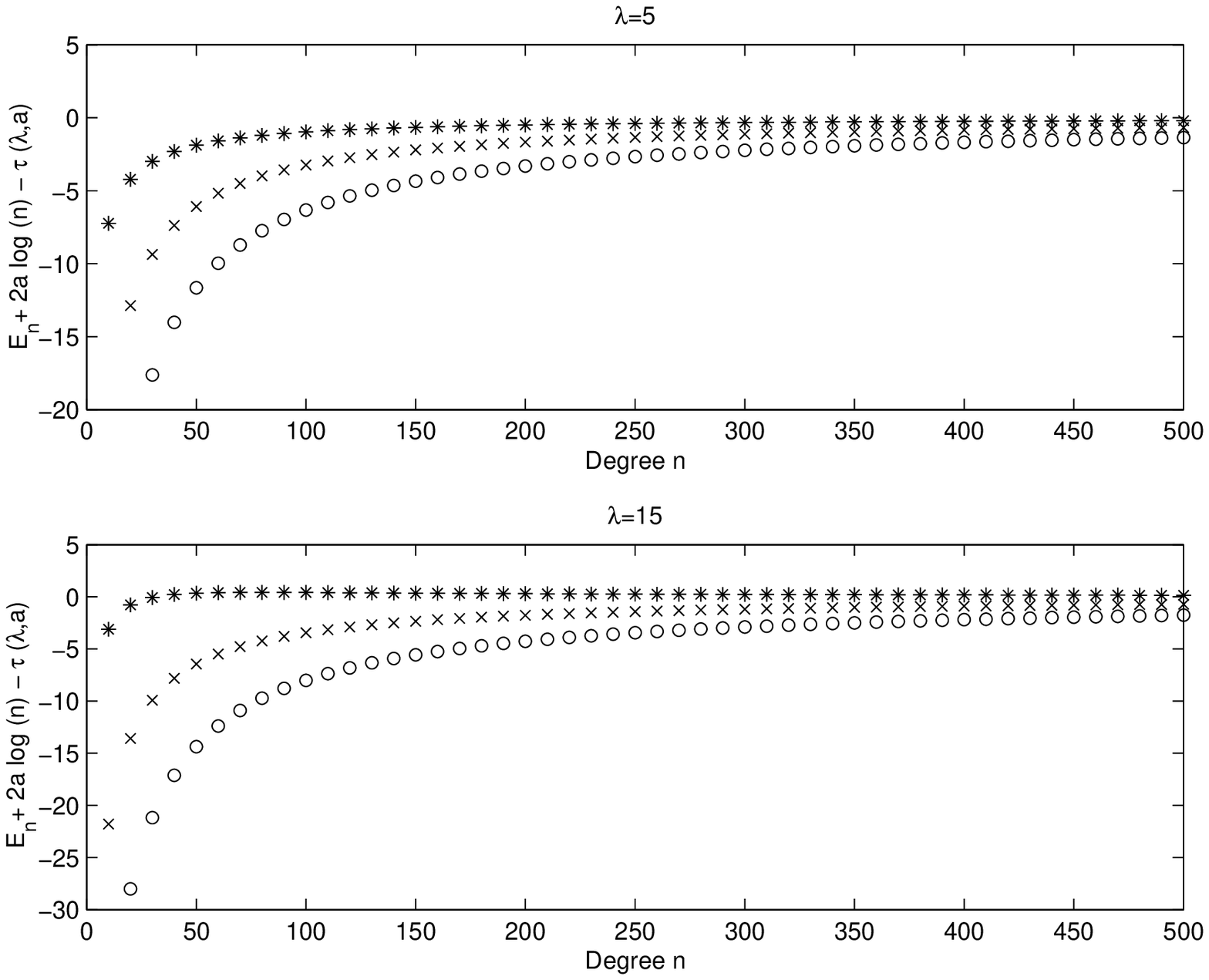}
\caption{$E_n(w_{\lambda}(\cdot;a))-\widetilde
E_n(w_{\lambda}(\cdot;a))$ for $n =1, 2, \dots, 500$, with $\lambda
=5$ (upper) and $\lambda=15$ (lower). We use the values $a=5$
(`\texttt{*}'), $a=10$ (`\texttt{x}') and $a=15$ (`\texttt{o}').
\label{grafPoll2}}
\end{figure}

For illustration, we have computed the entropy
$E_n(w_{\lambda}(\cdot;a))$ for $n =1, 2, \dots, 500$, $\lambda =5,
15$, and $a=5, 10, 15$ (Fig.\ \ref{grafPoll1}), using the numerical
algorithm from \cite{Buyarov04}, which admits as the only input data
the expression of the recurrence coefficients $a_n$ in
\eqref{CoefRecPoll}. For comparison, in Fig.\ \ref{grafPoll2} we
plot the difference $E_n(w_{\lambda}(\cdot;a))-
\widetilde{E}_n(w_{\lambda}(\cdot;a))$, where
$$
\widetilde{E}_n(w_{\lambda}(\cdot;a))\isdef -2a\log(n)+\tau (\lambda
;a)\,.
$$

In \cite{Levin03}, a rather general result about the leading term of
the asymptotics of $E_n$ has been established, that we state here
for the symmetric case: for even weights functions $w$ on $[-1,1]$
which belong to the class $\mathcal{F}(C^2+)$ introduced in
\cite{Levin01} (and whose definition we recall in Section
\ref{SecPeso}), if an additional assumption on the behavior of $w$
at $\pm 1$ (see Eq.\ (1) in \cite{Levin03}) is satisfied, then
\begin{equation}\label{LevinsAsymptotics}
E_n(w)=-\frac{2}{\pi}\int_{-\alpha_{n}}^{\alpha_n}
\frac{Q(x)}{\sqrt{\alpha_n^2-x^2 }}\, dx \, \big(1+o(1)\big), \quad
n\to \infty\,,
\end{equation}
where
\begin{equation}\label{defQ}
Q(x)\isdef -\frac{1}{2} \log\frac{w(x)}{w(0)}
\end{equation}
is the ``external field'', associated with the weight $w$, and
$\alpha _n$ is the \emph{Mahskar-Rakhmanov-Saff number} (or MRS
number), defined as the unique solution of the integral equation
\begin{equation}\label{defMRS1}
\frac{1}{\pi}\int_{-\alpha_{n}}^{\alpha_n}\frac{x~Q'(x)}{
\sqrt{\alpha_n^2-x^2}}~dx=n
\end{equation}
(see e.g.\ \cite{Saff:97} for details).
Unfortunately, the assumptions from \cite{Levin03} on the behavior
of $w$ at $\pm 1$ are not fulfilled by the symmetric Pollaczek
weights. Nevertheless, the result of Corollary \ref{TeoAsPoll} above
shows that the assertion in \cite{Levin03} is still valid:
\begin{cor}\label{CorAsPoll}
For the symmetric Pollaczek weight $w=w_{\lambda}(\cdot ;a)$, with
$a\geq 0$ and $\lambda\geq 1$, formula \eqref{LevinsAsymptotics}
holds.
\end{cor}
\begin{remark}
We should observe however that the result in \cite{Levin03} is
sharp: there the $o(1)$ term of \eqref{LevinsAsymptotics} has a
power decay, fact which is not true for the Pollaczek weight, where
the decay is logarithmic.
\end{remark}

Finally, asymptotic formula \eqref{asymptoticFormulaforE} has a
curious interpretation in terms of the behavior of the mutual energy
of two relevant sequences of probability measures on $[-1,1]$
associated with $p_n$'s:
 \begin{equation}\label{zerocounting}
\rho  _n =\frac{1}{n}\,  \sum_{j=1}^n \delta_{\zeta_j^{(n)}} \qquad
\text{and} \qquad d \nu_n(x)=p^2_n(x)w(x)\, d x\,,
\end{equation}
where $-1<\zeta_1^{(n)}<\dots < \zeta_n^{(n)}<1$ are the zeros of
the polynomial $p_n$. Both measures are standard objects of study in
the analytic theory of orthogonal polynomials. For instance, the
normalized zero counting measure $\rho _n$ is closely connected with
the $n$-th root asymptotics of $p_n$, while $\nu_n$ is associated
with the behavior of the ratio $p_{n+1}/p_n$ as $n \to \infty$ (see
\cite{Mate, Rakhmanov:77}).

If $\rho$ and $\nu$ are Borel (generally speaking, real signed)
measures on $\C$, we denote by
\begin{equation}\label{def_Log_potential}
 V^\rho(z) \isdef - \int \log |z-t|\,
d \rho (t)
\end{equation}
the logarithmic potential of $\rho$, and by
$$ I[\nu, \rho] \isdef \int
V^\nu(z) \, d\rho(z)= - \iint \log |z-t|\, d \nu (t) \, d\rho(z) \,,
$$ the mutual energy of $\nu$ and $\rho$. The latter is connected
with the entropy \eqref{def_E} by the formula
\begin{equation}\label{kullback2}
E_n =  - 2 \log \gamma _n +  2\sum_{j=1}^n V^{\nu_n} (\zeta_j^{(n)})
=- 2 \log \gamma _n + 2n \, I[\rho  _n, \nu _n]\,,
\end{equation}
where
\begin{equation}\label{leadingCoeff}
\gamma_n=2^n\left(\frac{(\lambda+a+1)_n~(\lambda+a)_n}{n!~(2\lambda)_n}\right)^{1/2}>0
\end{equation}
is the leading coefficient of $p_n$, and
$(z)_n=\Gamma(z+n)/\Gamma(z)$ denotes as usual the Pochhammer's
symbol. It is well known that as long as the orthogonality weight
$w>0$ a.e.\ on $[-1,1]$, both $\rho _n$ and $\nu_n$ tend (as $n \to
\infty$) in the weak-* sense to the Chebyshev (equilibrium)
distribution of the interval, which implies that $ \lim_{n\to
\infty} I[\rho  _n, \nu _n]=\log(2)$. In \cite{Beckermann04} a
rather surprising ``universal'' behavior of the next term of the
asymptotic expansion of $I[\rho  _n, \nu _n]$ was observed. Namely,
if the orthogonality weight satisfies the Szeg\H{o} condition
\eqref{szego} with an additional assumption on the growth of the
sequence of $\{ p_n\}$ on $[-1,1]$, then $ I[\rho _n, \nu _n] =
\log(2)-1/(2n)+o(1/n)$, $n \to \infty$. However, this result is no
longer valid for the Pollaczek polynomials, as it follows from
formula \eqref{asymptoticFormulaforE}:
\begin{cor}\label{CorEnergy}
For the symmetric Pollaczek weight $w=w_{\lambda}(\cdot ;a)$, with
$a\geq 0$ and $\lambda\geq 1$,
$$
I[\rho _n, \nu _n] =
\log(2)-\frac{1-2a}{2n}+o\left(\frac{1}{n}\right)\,, \quad n \to
\infty\,.
$$
\end{cor}
Observe that the second term of asymptotics is still independent of
the main parameter $\lambda $, and matches the result in
\cite{Beckermann04} for $a=0$.

The structure of this paper is as follows.  In Section \ref{SecPeso}
we gather some technical facts about the weight function
$w_{\lambda}(\cdot ;a)$; in particular, we show that this weight
does not satisfy \eqref{szego}, and for $\lambda \geq 1$, it belongs
to the class $\mathcal{F}(C^2+)$. Section \ref{sec:equilibrium}
contains some results about the equilibrium measure of total mass
$n$ in the external field $Q$; it is needed for the proof of Theorem
\ref{thFn} and Corollary \ref{TeoAsPoll} (Section \ref{SecFn}). We
defer the proof of the asymptotics of the sequence $\{G_n\}$
(Theorem \ref{thGn}) to Section \ref{SecGn}.  Finally, corollaries
\ref{CorAsPoll} and \ref{CorEnergy} are established in Section
\ref{SecCorLevin}.

\section{The weight function}\label{SecPeso}
As a first step in our analysis we study the behavior of the
symmetric Pollaczek weight function $w_{\lambda}(\cdot ;a)$ defined
in \eqref{WeightPollaczek}, which we denote simply by $w$ whenever
it cannot lead us into confusion. Using the notation
$$
t \isdef \frac{ax}{\sqrt{1-x^2}}\,, \quad x\in \Delta \isdef
[-1,1]\,,
$$
we rewrite its definition as
\begin{equation}\label{pesoPoll}
w(x)=\frac{2^{2\lambda}~(\lambda+a)}{2
\pi~\Gamma(2\lambda)}(1-x^2)^{\lambda-1/2}~e^{(2\arccos
x-\pi)t}~\left|\Gamma(\lambda+i t)\right|^2, \quad x \in \Delta\,.
\end{equation}
This is an even function on $\Delta$, strictly positive in $(-1,1)$,
but vanishing at the end points. The fast (exponential) decay at
$\pm 1$ is precisely the reason why $w$ does not satisfy
\eqref{szego}.
Indeed, using the asymptotic formula (6.1.40) from
\cite{abramowitz/stegun:1972} (see also \cite[\S 2.11]{Luke75}), we
can easily obtain that
\begin{equation*}
2\log\left|\Gamma\left(\lambda+ti\right)\right|= -\pi t+
2(\lambda-1/2)\log
t+2\log\sqrt{2\pi}+\frac{\lambda(\lambda-1)(2\lambda-1)}{6t^2}+\mathcal
O(t^{-4}),\end{equation*} when $x\rightarrow 1^-$ (which denotes in
what follows the one-sided limit from the left). Thus,
(\ref{pesoPoll}) yields that
\begin{align}
\log w(x)&= \log\left(\frac{2^{2\lambda}~(\lambda+a)}{2
\pi~\Gamma(2\lambda)}\right)+(\lambda-1/2)\log(1-x^2)+(2\arccos
x-\pi)t\nonumber\\&\quad+\log\left|\Gamma(\lambda+i t)\right|^2\nonumber\\
&=-2\pi t+(\lambda-1/2)\log(x^2)+2t \arccos x
+\log\left(\frac{2^{2\lambda}~(\lambda+a)a^{2\lambda-1}}{2
\pi~\Gamma(2\lambda)}\right)\nonumber\\& \quad+2\log\sqrt{2\pi}
+\frac{\lambda(\lambda-1)(2\lambda-1)}{6t^2}+\mathcal
O(t^{-4}),\quad x \rightarrow 1^-. \label{AsLogPesoPoll}
\end{align}
Since
\begin{equation}\label{arcos}
2t \arccos x=2a-\frac{2a}{3}(1-x^2)+\mathcal O(1-x^2)^2,\qquad
x\rightarrow 1^-,
\end{equation}
we obtain that
$$\log w(x)=-2\pi t+\mathcal O(1),\qquad x\rightarrow
1^-.
$$
In consequence,
\begin{equation}\label{AsPesoPoll} w(x)=\exp(-2\pi
|t|+\mathcal O(1))=\exp\left(\frac{-2\pi a
|x|}{\sqrt{1-x^2}}+\mathcal O(1)\right), \quad |x|\to 1^-\,,
\end{equation}
showing that for this weight the integral in \eqref{szego} is
divergent.
The previous analysis motivates the introduction of functions $w_0$
and $s$ on $\Delta$, such that
\begin{equation}\label{defw0}
w_0(x) \isdef e^{-2\pi|t|}=\exp\left(-~\frac{2\pi a
|x|}{\sqrt{1-x^2}}\right), \quad \text{and} \quad
w(x)=w_0(x)e^{s(x)}.
 \end{equation}
\begin{lm}\label{prolims}
Function $s \in C^\infty(-1,1)$ is even and continuous in $[-1,1]$.
\end{lm}
\begin{proof} We need to check only the existence of finite limits of $s$ at $\pm 1$ (the rest
is trivial). From (\ref{AsLogPesoPoll}) it is clear that
\begin{align*}
s(x) =\log\frac{w(x)}{w_0(x)} =&
\log\left(\frac{2^{2\lambda}~(\lambda+a)}{2
\pi~\Gamma(2\lambda)}~a^{2\lambda-1}\right)+(\lambda-1/2)\log
(x^2)\\& +2 t\,\arccos
x+\frac{\lambda(\lambda-1)(2\lambda-1)}{6t^2}+\mathcal O(t^{-4}),
\end{align*}
and using (\ref{arcos}) we get
\begin{align*}
s(x)= \log\left(\frac{2^{2\lambda}
(\lambda+a)}{\Gamma(2\lambda)}~a^{2\lambda-1}e^{2a}\right)+\mathcal
O(x-1), \quad x\to 1^-,
\end{align*}
which concludes the proof.
\end{proof}

\begin{prop}
\label{prop:weight} For $\lambda\geq 1$, the weight
$w=w_{\lambda}(\cdot ;a)$ belongs to the class $\mathcal{F}(C^2+)$.
\end{prop}
Recall that $ w \in \mathcal{F}(C^2+)$ (see \cite{Levin01}) if the
corresponding external field $Q$ introduced in \eqref{defQ} is
positive and verifies the following conditions:
\begin{enumerate}
\item[a)] $Q'$ is continuous in $\Delta$.
\item[b)] $Q'{'}$ exists and is positive in $\Delta \setminus
\{0\}$.
\item[c)]$\lim_{|x|\to 1^-}Q(x)=\infty.$
\item[d)] Function $$T(x)\isdef \frac{xQ'(x)}{Q(x)},$$ is
cuasi-increasing in $(0,1)$ cuasi-decreasing in $(-1,0)$, and
\begin{equation}\label{boundForT}
T(x)\geq\Lambda>1, \quad x\in (-1,1).
\end{equation}
\item[e)] There exists $C_1>0$ such that
$$\frac{Q'{'}(x)Q(x)}{(Q'(x))^2}\leq C_1, \quad x\in (-1,1).$$
\item[f)] There exist a compact subinterval $J$, contained in $(-1,1)$,
and $C_2>0$ such that
  $$\frac{Q'{'}(x)Q(x)}{(Q'(x))^2}\geq C_2,$$
for all $x\in (-1,1)\setminus J$, except a subset with zero measure.
\end{enumerate}
A function $f:I \rightarrow [0,+\infty)$ is cuasi-increasing if
$\forall x<y\in I,~ \exists C>0$, such that $f(x)<C~f(y).$

We prove Proposition \ref{prop:weight} in several steps.
\begin{lm}
\label{lemma:integral} If  a function $f\in C[0,+\infty)$ is
positive and decreasing, then for $t>0$,
\begin{equation}\label{IntPos}
I(t)\isdef \int_0^{+\infty} f(u) \sin (ut)\, du \in ( 0, +\infty]\,.
\end{equation}
 \end{lm}
\begin{proof}  With
the change of variable $v=ut$ we can write
 $$I(t)=\frac 1t\int_0^\infty f(v/t)~\sin (v)~ dv.$$
We denote $g(v)=f(v/t)$; then
\begin{align*}
t I(t)& =\sum_{k=0}^\infty\int_{k\pi}^{(k+1)\pi}g(v)~\sin v~dv=
 \sum_{k=0}^\infty\int_{0}^{\pi}g(k\pi+v)~\sin(k\pi+v)~dv\nonumber\\
 &=\sum_{k=0}^\infty\int_{0}^{\pi}g(2k\pi+v)\sin(2k\pi+v)dv \\ &\qquad +
 \sum_{k=0}^\infty\int_{0}^{\pi}g((2k+1)\pi+v)\sin((2k+1)\pi+v)dv\nonumber\\
 &=\sum_{k=0}^\infty\int_{0}^{\pi}g(2k\pi+v)~\sin(v)~dv-
 \sum_{k=0}^\infty\int_{0}^{\pi}g((2k+1)\pi+v)~\sin(v)~dv\nonumber\\
 &=\sum_{k=0}^\infty\int_{0}^{\pi}\bigg(g(2k\pi+v)-g((2k+1)\pi+v)\bigg)\, \sin(v)~dv.
\end{align*}
Since each integral in the series is strictly positive, it proves
(\ref{IntPos}).
\end{proof}
Now we gather some properties of the digamma and trigamma functions
in the following technical lemma:
\begin{lm}\label{proPhi}
For the digamma function $\psi(x)=\Gamma'(x)/\Gamma(x)$  the
following statements hold: for $\lambda\geq 1$,
\begin{enumerate}
\item[i)] $\Re\psi'(\lambda+it)$ is a strictly positive even function of $t\in \R$, and
\begin{equation}\label{limitRePsi}
\lim_{t\rightarrow \pm\infty}t \Re\psi'(\lambda+it)=0.
\end{equation}
\item[ii)]
$\Im\psi(\lambda+it)$ is an odd function of $t\in \R$, strictly
positive in $(0,+\infty)$, and
\begin{equation}\label{limitImPsi}
\Im\psi(\lambda\pm it)=\pm\frac{\pi}{2}+\mathcal O(t^{-1})\,, \quad
t\rightarrow +\infty.
\end{equation}
\end{enumerate}
\end{lm}
\begin{remark}
Limits \eqref{limitRePsi}--\eqref{limitImPsi} are valid in fact for
$\lambda
>0$.
\end{remark}
\begin{proof} The symmetry of both the real and the imaginary parts
of $\psi(\lambda +i t)$ follows from the well known property
$$
\psi(\overline{z})=\overline{\psi(z)}\,;
$$
so, we restrict our attention to $t>0$.

For \emph{i)} we  consider the integral representation of the
trigamma function,
\begin{equation}\label{intgralRepresentationPsi}
\psi'(z) =\int_0^{+\infty} e^{-uz}\frac{u}{1-e^{-u}}du,
\end{equation}
(see e.g.\ \cite[formula 3.41.371.6]{Gradshtein95}), from where
$$
\Re\psi'(\lambda+it)=\int_0^\infty e^{-u\lambda}\frac{u}{1-e^{-u}} \cos (ut)~du.
$$
Integrating by parts it can be reduced to
$$
\Re\psi'(\lambda+it)=\frac{1}{t}\int_0^\infty e^{-\lambda
u}~\frac{ue^{-u}-(1-\lambda
 u)(1-e^{-u})}{(1-e^{-u})^2}\, \sin (ut)~du.
$$
It is easy to check that for $\lambda\geq 1$, function
$$f(u)=e^{-\lambda u}~\frac{ue^{-u}-(1-\lambda
 u)(1-e^{-u})}{(1-e^{-u})^2},$$
is positive and decreasing on $(0,+\infty)$. Hence, the first part
of \emph{i)} follows from Lemma \ref{lemma:integral}.
On the other hand, by \cite[formula 6.4.12]{abramowitz/stegun:1972},
$$
\psi'(z)=\frac{1}{z} +\mathcal O \left(\frac{1}{z^2}\right)\,, \quad
z \to \infty\,,\quad |\arg z |<\pi\,,
$$
so that
$$
(\lambda + it)\psi'(\lambda + it)=1+\mathcal O
\left(\frac{1}{t}\right)\,, \quad t \to \pm \infty\,.
$$
In particular,
\begin{align*}
 0& =\lim_{t\to \pm \infty} \Im \left((\lambda + it)\psi'(\lambda +
it) \right)=\lambda \lim_{t\to \pm \infty} \Im \left(\psi'(\lambda +
it) \right)+ \lim_{t\rightarrow \pm\infty}t \Re\psi'(\lambda+it)\\
&=\lim_{t\rightarrow \pm\infty}t \Re\psi'(\lambda+it)\,,
\end{align*}
which proves \eqref{limitRePsi}.
\medskip
For \emph{ii)} we may use the series expansion \cite[formula
6.3.16]{abramowitz/stegun:1972},
$$
\psi(1+z)=-\gamma+\sum_{k=1}^\infty \frac{z}{k(k+z)}\,, \quad
-z\notin \N\,,
$$
according to which
$$
\Im \psi(1+z)=-\sum_{k=1}^\infty \Im \left(
\frac{1}{k+z}\right)=\sum_{k=1}^\infty \frac{\Im(z)}{|k+z|^2} >0
\quad \text{if } \Im(z)>0\,.
$$
%
Finally, by the asymptotic formula \cite[formula
6.3.18]{abramowitz/stegun:1972},
$$
\psi(z)=\log z+\mathcal O \left(\frac{1}{z}\right)\,, \quad z \to
\infty\,,\quad |\arg z |<\pi\,,
$$
from which \eqref{limitImPsi} is immediate.
\end{proof}
Now we are ready to analyze whether the weight belongs to $\mathcal
F(C^2+)$. Conditions a) and b) are a straightforward consequence of
the following statement:
\begin{lm}
\label{lemma:Q} The even function $Q\in C^\infty(-1,1)$ satisfies
$$
\frac{d^k \, Q(x)}{dx^k}  >0 \quad \text{for  $ x\in (0,1)$ and
$k=0,1,2$.}
$$
\end{lm}
\begin{proof}
By definition, $Q(0)=0$, and by symmetry, $Q'(0)=0$, so
$$
Q''(x)>0 \quad \Rightarrow \quad Q'(x)>0 \quad \Rightarrow \quad
Q(x)>0\,, \quad x\in (0,1)\,.
$$
But for $x\in (0,1)$ we have $t>0$, and
\begin{equation}\label{QdosPrime}
\begin{split}
Q''(x) &=\frac{(\lambda-1/2)(1+x^2)}{(1-x^2)^2}
 +\frac{ax^2+2a}{(1-x^2)^2}+\frac{3ax}{(1-x^2)^{5/2}}~\frac{1}{2}(\pi-2\arccos x) \\
&\quad+\frac{3ax}{(1-x^2)^{5/2}}\Im\psi(\lambda+it)+
\frac{a^2}{(1-x^2)^3} \Re\psi'(\lambda+it)>0\,,
\end{split}
\end{equation}
where we have used Lemma \ref{proPhi}.
\end{proof}
Since by (\ref{AsPesoPoll}),
\begin{equation}\label{AsQ}
Q(x)=-~\frac{1}{2}\log w(x)+\frac{1}{2}\log
w(0)=\frac{1}{2}~\frac{2\pi a |x|}{\sqrt{1-x^2}}+\mathcal O(1),
\quad |x|\rightarrow 1^-,
\end{equation}
condition c) also trivially holds.

We turn now to the even function
$$
T(x)=\frac{x Q'(x)}{Q(x)};
$$
let us show that it is cuasi-increasing in $(0,1)$. Since $T$ is
continuous and positive on the bounded interval $(0,1)$, it is
sufficient to show that it does not blow up at the left end point,
nor it vanishes at the right one. Recall that
$$
Q(0)=Q'(0)=0,\quad Q'{'}(0)=( \lambda-1/2)+2a+a^2\psi'(\lambda)>
0\,;
$$
in particular,
\begin{equation}\label{localQ}
    Q(x)=\frac{( \lambda-1/2)+2a+a^2\psi'(\lambda)}{2}\,x^2
    +\mathcal O(x^3)\,, \quad x\to 0\,.
\end{equation}
Hence,
\begin{equation} \label{limTen0}
 \lim_{x\rightarrow
0}T(x)=\lim_{x\rightarrow
 0} \frac{x  Q'(x)}{Q(x)}=\lim_{x\rightarrow
 0}
 \frac{Q'(x)/x}{Q(x)/x^2}=\frac{Q'{'}(0)}{\frac{1}{2}Q'{'}(0)}=2>0\,.
\end{equation}
On the other hand,
\begin{equation}\label{derQ}
Q'(x)=\frac{(\lambda-1/2)x}{1-x^2}+\frac{ax}{1-x^2}+\frac{1}{2}~
\frac{(\pi-2\arccos
x)~a}{(1-x^2)^{3/2}}+\frac{a}{(1-x^2)^{3/2}}\Im\psi(\lambda+i t),
\end{equation}
and by Lemma \ref{proPhi},
\begin{equation}\label{asderQ}
\lim_{x\rightarrow 1^-}(1-x^2)^{3/2} ~Q'(x)=a\pi\,.
\end{equation}
Together with (\ref{AsQ}) it shows that
\begin{equation}\label{divT}
\lim_{x\rightarrow 1^-}T(x)=\lim_{x\rightarrow 1^-}\frac{x
Q'(x)}{Q(x)} =\lim_{x\rightarrow 1^-}\frac{1}{1-x^2}~\frac{\pi
a}{\pi a}=+\infty.
\end{equation}
In conclusion, $T$ is cuasi-increasing in $(0,1)$.

On the other hand, if $\zeta\in (0,1)$ is a local minimum of $T$,
then
$$
\frac{d \log(T(x)}{dx}\bigg|_{x=\zeta}=0\quad \Rightarrow \quad
\frac{1}{\zeta }+\frac{Q''}{Q'}(\zeta )-\frac{Q'}{Q}(\zeta )=0\,,
$$
or equivalently,
$$
T(\zeta )=1+\frac{\zeta Q''(\zeta )}{Q'(\zeta )}>1\,.
$$
Taking into account also the behavior at $x=0$ and $x=1$ (see
\eqref{limTen0} and \eqref{divT}), we obtain \eqref{boundForT}.

Finally, let us check conditions e) and f). Denote
$$
H(x)\isdef \frac{Q(x) Q''(x)}{(Q'(x))^2}\,.
$$
This is an even, continuous and positive function on $(-1,1)$, with
$ H(0)=1/2$, where we have used \eqref{localQ}.
On the other hand,
$$\lim_{x\rightarrow 1^-} (1-x^2)^{1/2} Q(x) =\lim_{x\rightarrow
1^-}(1-x^2)^{3/2} Q'(x)=a\pi\,,
$$
where we have used (\ref{AsQ}) and (\ref{asderQ}). Also from
\eqref{QdosPrime} it follows that
$$\lim_{x\rightarrow
1^-}(1-x^2)^{5/2} Q''(x)= \frac{3}{2}\, a\pi+\frac{3}{2}\, a
\pi+a\lim_{x\rightarrow 1^-}t\Re\psi'(\lambda+it),
$$
and by \eqref{limitRePsi}, $\lim_{x\rightarrow
1^-}(1-x^2)^{5/2}~Q'{'}(x)=3 a\pi$. Gathering these identities we
obtain that
$$ \lim_{x\rightarrow 1^-}H(x)=\lim_{x\rightarrow 1^-} \frac{(1-x^2)^{1/2} Q(x)\,
(1-x^2)^{5/2} Q''(x)}{ (1-x^2)^{3 } (Q'(x))^2}=3\,.$$ Hence, $H$ can
be extended as a strictly positive and continuous (and thus,
uniformly continuous) function on $[-1,1]$; from this fact
conditions e) and f) follow automatically. This concludes the proof
of Proposition \ref{prop:weight}.

\section{Equilibrium measure}\label{sec:equilibrium}

The equilibrium measure $\mu_n$ on $[-1,1]$ of total weight $n$ in
the external field $Q$ plays a prominent role in the asymptotics of
the orthogonal polynomials, and we gather in this section some of
its properties needed further. By Lemma \ref{lemma:Q}, function
$Q\in C^\infty(-1,1)$ is strictly convex, and
$Q(-1^+)=Q(1^-)=+\infty$. In consequence, (see e.g. \cite{Saff:97}),
$\mu_n$ is absolutely continuous and supported on the interval
$[-\alpha _n, \alpha _n]$, where $\alpha _n$ is the MRS number,
defined by \eqref{defMRS1}. If we denote by $\sigma _n(x)$ the
density ($\mu_n'$) of $\mu_n$, then
$$
\int_{-\alpha _n}^{\alpha _n} \sigma_n(x)\, dx=n\,,
$$
and the characterizing property of the equilibrium is
$$
V^{\mu _n}(x)+Q(x) \begin{cases} = b_n (=\const), & x\in  [-\alpha
_n, \alpha _n], \\ > b_n, & \alpha _n<|x|\leq 1\,,
\end{cases}
$$
where $ V^{\mu_n}$ is the logarithmic potential of $\mu_n$ (cf.\
\eqref{def_Log_potential}). We analyze first the asymptotic behavior
of $\{\alpha _n\}$, needed in the proof of Corollary
\ref{CorAsPoll}. It is known that $\alpha _n \to 1^-$, and even
more, that $1-\alpha _n=\mathcal O(1/n)$ (see \cite[\S
1.6]{Levin01}), but we are looking for a more precise information.
The following technical lemma is useful for the estimation of the
behavior of the integral in \eqref{defMRS1}:
\begin{lm}\label{proIntElipticas}
Let $f(u,x)$ be defined for $u,x\in[0,1]$ with $f$ and $
\partial f/\partial x$ continuous in $[0,1]^2$, $f(1,1)\neq 0$; then
when $u\rightarrow 1^-$,
\begin{align} \label{ellInt1}
\int_0^1\frac{f(u,x)}{\sqrt{1-u^2x^2}}\, \frac{dx}{\sqrt{1-x^2}} & =
-\frac{f(1,1)}{2}\log(1-u)\, \big( 1+o(1)\big),\\\label{ellInt2}
\int_0^1\frac{f(u,x)}{(1-u^2x^2)}~\frac{dx}{\sqrt{1-x^2}}&=
\frac{\pi f(1,1)}{2\sqrt{2}}~\frac{1}{\sqrt{1-u}}\, \big(
1+o(1)\big),\\\label{ellInt3}
 \int_0^1\frac{f(u,x)}{(1-u^2x^2)^{3/2}}~\frac{dx}{\sqrt{1-x^2}}&=
 \frac{f(1,1)}{2}~\frac{1}{1-u}\, \big( 1+o(1)\big).
\end{align}
\end{lm}
Formula \eqref{ellInt1} appears in \cite[formula
8.113.3]{Gradshtein95} for $f(u,x)=1$; the proof is standard, and we
omit it here for the sake of brevity.

Now we can obtain the first two terms of the asymptotics of $\alpha
_n$:
\begin{prop}\label{proMRS}
The MRS numbers  $\alpha_n$  satisfy
$$ \alpha_n=1- \frac{a}{n}+o\left( \frac{1}{n}\right)\,, \quad
n\rightarrow\infty\,.$$
\end{prop}
\begin{proof}
Formula \eqref{defMRS1} defining the MRS numbers may be rewritten as
\begin{equation}\label{defMRS2}
\frac{2}{\pi}~\int_{0}^1\frac{\alpha_nx~Q'(\alpha_n
x)}{\sqrt{1-x^2}}~dx=n\,,
\end{equation}
which motivates the study of the asymptotics (as  $u\rightarrow
1^-$) of the integral
\begin{align*}
\frac{2}{\pi}~&\int_0^1\frac{ux~Q'(ux)}{\sqrt{1-x^2}}~dx\\
 &=\frac{1}{\pi}~\int_0^1\frac{ux~\left((\lambda-1/2)2ux+2aux\right)}{(1-u^2x^2)~\sqrt{1-x^2}}~dx
 +\frac{1}{\pi}~\int_0^1\frac{ux~(-2)a~\arccos(ux)}{(1-u^2x^2)^{3/2}~\sqrt{1-x^2}}~dx\\
 &\quad+\frac{1}{\pi}~\int_0^1\frac{ux\left(\pi a +2a\Im \psi\left(\lambda+i \frac{aux}{\sqrt{1-u^2x^2}}\right)\right)}
 {(1-u^2x^2)^{3/2}~\sqrt{1-x^2}}~dx=I_1(u)+I_2(u)+I_3(u).
\end{align*}
Using Lemma \ref{proIntElipticas} we have that for $u\to 1^-$,
\begin{align*} I_1(u)&=
\frac{1}{\pi}((2\lambda-1)+2a)\frac{\pi}{2\sqrt{2}}\frac{1}{\sqrt{1-u}}\,\big(1+o(1) \big)
=o\left(\frac{1}{1-u}\right), \\
I_3(u)& = \frac{1}{\pi}~(\pi a+2a\pi/2)\frac{1}{2}~\frac{1}{1-u}\,
\big(1+o(1) \big) =\frac{a}{1-u}\, \big(1+o(1) \big).
\end{align*}
Moreover, since
$$f(u,x)=\frac{\arccos(ux)}{\sqrt{1-u^2x^2}}\, ,$$ satisfies the conditions of
Lemma \ref{proIntElipticas}, we have that
\begin{equation*}\label{I-MRS2}
I_2(u)= -\frac{1}{\pi}\int_0^1\frac{u x 2 a
f(u,x)}{(1-u^2x^2)\sqrt{1-x^2}}=o\left(\frac{1}{1-u}\right), \quad
u\to 1^-.
\end{equation*}
Summarizing,
$$\frac{2}{\pi}~\int_0^1\frac{ux~Q'(ux)}{\sqrt{1-x^2}}~dx = \frac{a}{1-u}\,
\big(1+o(1) \big),$$ and equation \eqref{defMRS2} for the MRS
numbers can be rewritten as
 $$
 \frac{a}{1-\alpha_n}\,
\big(1+o(1) \big)=n,$$
 which proves the statement.
\end{proof}
\medskip
We turn now to the analysis of the density $\sigma_n$. It is
convenient to introduce the normalized translation of $\sigma_n$ to
$\Delta$,
\begin{equation}\label{defSigma_norm}
\sigma_n^*(u)=\frac{\alpha _n}{n}\, \sigma _n(\alpha _n u)\,, \quad
u \in \Delta\,,
\end{equation}
as well as the cumulative distribution
\begin{equation}\label{defPhi}
\Phi_n(\theta) \isdef \pi\int_{\cos\theta}^1\sigma_n^*(t)dt\,,
\end{equation}
which is obviously a smooth and strictly increasing function on
$[0,\pi]$; moreover, $\Phi_n:\, [0,\pi] \to [0,\pi]$ is a bijection,
and the inverse function $\Phi_n^{[-1]}$ exists (observe that our
definition differs in normalization from that used in \cite{Levin01,
Levin03}). We summarize some properties of $\Phi_n$ in the following
lemma.
\begin{lm}\label{lmPhin} For $\Phi_n$ defined in \eqref{defPhi},
\begin{itemize}
\item[i)]$\Phi_n'(\theta) \to 1 $ pointwise in $(0,\pi)$.
\item[ii)]$$\int_0^\pi\left|\frac{1}{\Phi_n'(\Phi_n^{[-1]}(\eta))}-1\right|~d\eta=o(1),\qquad
n\to\infty.$$
 \item[iii)] $\Phi_n(\theta) \to \theta$ and
$\Phi_n^{[-1]}(\theta) \to \theta$ as $n\to \infty$ uniformly in
$[0,\pi]$.

 \end{itemize}
\end{lm}
\begin{proof}
 Let $0< \varepsilon <1/2$, and $x=\cos
 \theta\in (-1+2 \varepsilon ,
1-2 \varepsilon )$. By Lemma 6.5 of \cite{Levin01},
\begin{equation*}\label{reprePhi'New}
\Phi_n'(\theta)-1  =\pi\sigma_n^*(\cos\theta)\sin\theta-1=
\frac{\alpha _n}{\pi n }\, PV\, \int_{-1}^1
 \frac{Q'(\alpha_n u) \sqrt{1-u^2}}{u-x}\, du  \,,
\end{equation*}
where $PV$ means the principal value of integral. Hence,
\begin{equation}\label{reprePhi'}
\begin{split}
\Phi_n'(\theta)-1& =
 \frac{\alpha _n}{\pi n }\, \left(
 \int_{-1+\varepsilon}^{1-\varepsilon}
 \frac{Q'(\alpha_n u) \sqrt{1-u^2}-Q'(\alpha_n x) \sqrt{1-x^2}}{u-x}\, du  \right. \\
 &+  Q'(\alpha_n x) \sqrt{1-x^2} PV\, \int_{-1+\varepsilon}^{1-\varepsilon}
 \frac{du}{u-x}+ \int_{-1 }^{-1+\varepsilon}
 \frac{Q'(\alpha_n u) \sqrt{1-u^2}}{u-x}\, du \\
 &\left. + \int_{1-\varepsilon }^{1}
 \frac{Q'(\alpha_n u) \sqrt{1-u^2}}{u-x}\, du\right) \,.
\end{split}
\end{equation}
The first two terms within parentheses in the right hand side of
\eqref{reprePhi'} are uniformly bounded. Let us estimate
$$
\left| \int_{1-\varepsilon }^{1} \frac{Q'(\alpha_n u)
\sqrt{1-u^2}}{u-x}\, du \right| \leq \frac{1}{\varepsilon }\,
\int_{1-\varepsilon }^{1}  Q'(\alpha_n u) \sqrt{1-u^2} \, du
$$
(the remaining term is analyzed in a similar fashion). Integrating
by parts,
$$
\int_{1-\varepsilon }^{1}  Q'(\alpha_n u) \sqrt{1-u^2} \,
du=-\frac{1}{\alpha _n}\, Q(\alpha_n (1-\varepsilon))
\sqrt{\varepsilon(2-\varepsilon)}+ \frac{1}{\alpha
_n}\,\int_{1-\varepsilon }^{1} \frac{u Q(\alpha_n u) }{\sqrt{1-u^2}}
\, du\,,
$$
and using \eqref{ellInt1} we get that
$$
\left|\Phi_n'(\theta)-1\right|=\mathcal O\left( \frac{\log
(n)}{n}\right)\,,
$$
which proves \emph{i)}.

On the other hand, by \cite[lemma 4.2 a)]{Levin03}, the sequence
$|\Phi_n'(\theta)|$ is uniformly bounded on $[0,\pi]$. Thus, by the
dominated convergence theorem,
\begin{align*}
\int_0^\pi\left|\frac{1}{\Phi_n'(\Phi_n^{[-1]}(\eta))}-1\right|~d\eta
& =
\int_0^\pi\left|\frac{1}{\Phi_n'(\theta)}-1\right|~\Phi_n'(\theta)d\theta
\\ & = \int_0^\pi\left|1-\Phi_n'(\theta)\right|~d\theta \rightarrow 0\,,
\quad n\to \infty\,.
\end{align*}

Furthermore, given $\theta\in [0,\pi]$,
 $$\left|\Phi_n(\theta)-\theta\right|=\left|\int_0^\theta(\Phi_n'(\eta)-1)d\eta\right|\leq
 \int_0^{\pi}\left|\Phi_n'(\eta)-1\right|d\eta=o(1),$$
and so the uniform convergence of $\Phi_n$ on $[0,\pi]$ follows
again  by the dominated convergence theorem.
Finally, if $\theta=\Phi_n(\eta)$, 
\begin{align*}
&\left|\Phi_n^{[-1]}(\theta)-\theta\right|=\left|\Phi_n^{[-1]}(\Phi_n(\eta))-\Phi_n(\eta)\right|=
\left|\eta-\Phi_n(\eta)\right|,
\end{align*}
showing that $\Phi_n^{[-1]}(\theta)$ converges uniformly to $\theta$
on $[0,\pi]$. This concludes the proof.
\end{proof}

\section{Asymptotics of $F_n$: proof of Theorem \ref{thFn}}\label{SecFn}

We follow the scheme of proof of \cite{Beckermann04}. We have
established already that the weight $w\in \mathcal F(C^2+)$;
furthermore, $w(x)>0$ on $(-1,1)$, and in consequence, it is an
Erd\H{o}s-Tur\'{a}n weight. One of the most relevant facts about these
weights is that the sequence $p_n^2(x)w(x)dx$ converges in the
weak-$*$ topology to the equilibrium (Robin) measure $\mu$ of
$[-1,1]$ (see \cite{Mate, Rakhmanov:77}). In other words, for any
$f\in C[-1,1]$,
\begin{equation}\label{convErdos}
\lim_{n\rightarrow\infty}\int_{-1}^1 f(x) p_n^2(x) w(x)~ dx=
\frac{1}{\pi}\, \int_{-1}^1 f(x)\frac{dx}{ \sqrt{1-x^2}}.
\end{equation}
We make use also of the  following technical lemmas. In the sequel
we write that $x_n \sim y_n$ if the ratios $x_n/y_n$ and $y_n/x_n$
are uniformly bounded in $n$.
\begin{lm}\label{lemaFnMRS}
\begin{equation*}
\int_{\Delta\setminus[-\alpha_n,\alpha_n]}\log(p_n^2(x)~w(x))~p_n^2(x)~w(x)~dx=o(1)\,,
\quad n \to \infty\,.
\end{equation*}
\end{lm}
\begin{proof}
Since $w\in \mathcal F(C^2+)$, by  \cite[Theorem 1.18]{Levin01},
$$
\sup_{x\in\Delta}\left|p_n(x)~\sqrt{w(x)}\right|\sim
n^{1/6}~(\alpha_n)^{-1/3}~ \left(\frac{T(\alpha_n)}{\alpha
_n}\right)^{1/6}\sim n^{1/3},
$$
where we have taken into account Proposition \ref{proMRS}. Thus, for
$\varepsilon>0$, there exists $C_1>0$ such that
$$\left(p_n^2(x)~w(x)\right)^{1+\varepsilon}\leq C_1\, n^{(2+2\varepsilon)/3}$$
and so, if $\varepsilon<1/2$, then
\begin{align*}
&\int_{\Delta\setminus[-\alpha_n,\alpha_n]}\log(p_n^2(x)~w(x))~p_n^2(x)~w(x)~dx\leq
\int_{\Delta\setminus[-\alpha_n,\alpha_n]}(p_n^2(x)~w(x))^{1+\varepsilon}~dx\\
\leq
&C_1~n^{(2+2\varepsilon)/3}~\int_{\Delta\setminus[-\alpha_n,\alpha_n]}dx
\leq
C_2~n^{\frac{2+2\varepsilon}{3}-1}=C_2~n^{\frac{-1+2\varepsilon}{3}}=o(1).
\end{align*}
For a lower bound it is sufficient to take into account that
function
\begin{equation}\label{defR}
\mathcal{R}(y)=y^2 \log(y^2)
\end{equation}
is bounded from below on $[0,+\infty)$, and thus
$$
\int_{\Delta\setminus[-\alpha_n,\alpha_n]}\log(p_n^2(x)~w(x))~p_n^2(x)~w(x)~dx\geq
C\, \int_{\Delta\setminus[-\alpha_n,\alpha_n]} dx=o(1) \,, \quad
n\to \infty\,,$$ which concludes the proof.
\end{proof}

\medskip

\begin{lm}\label{lemaFn}
\begin{equation}\label{limitPotentials}
\lim_{n\rightarrow\infty}
\int_{-\alpha_n}^{\alpha_n}\log(\sqrt{\alpha_n^2-x^2})~p_n^2(x)~w(x)~dx=-\log
(2).
\end{equation}
\end{lm}

\begin{proof}
Let
$$
\ell_n(x)\isdef \begin{cases} \log(\sqrt{\alpha_n^2-x^2}), & x\in
(-\alpha _n, \alpha _n), \\
0, & x\in [-1,1]\setminus (-\alpha _n, \alpha _n)\,.
\end{cases}
$$
Obviously, $\ell_n(x)\to\log(\sqrt{1-x^2})$ pointwise for $x\in
(-1,1)$. Furthermore, there exists $C_1>0$ such that for $t>1$ it
holds that $(\log(t))^3<C_1 t$. Then, by Lebesgue dominated
convergence theorem,
$$
\left\|\ell_n(x)- \log(\sqrt{1-x^2})\right\|_{L^3}=o(1)\,, \quad
n\to\infty
$$
(here and in the sequel $\|\cdot\|_{L^p}$ denotes the $p$-norm with
respect to the Lebesgue measure on $\Delta$).

Furthermore, from \cite[Theorem 13.6]{Levin01} it follows that the
sequence $\|p_n\sqrt{w}\|_{L^p}$ is uniformly bounded as long as
$p<4$ (and in particular, for $p=3$). Thus, by H\"{o}lder inequality,
\begin{equation*}
\begin{split}
 \left| \int_{-1}^{1}\left( \ell_n(x)- \log(\sqrt{1-x^2})\right)\,
p_n^2(x)~w(x)   dx\right| \\   \leq \left\|\ell_n(x)-
\log(\sqrt{1-x^2})\right\|_{L^3}\, \left\|p_n^2 ~w
\right\|_{L^{3/2}} \,,
\end{split}
\end{equation*}
which implies that
\begin{equation}\label{limitFirstIntegral}
\begin{split}
    \int_{-\alpha_n}^{\alpha_n}\log(\sqrt{\alpha_n^2-x^2})~p_n^2(x)~w(x)~dx
\\ =\int_{-1}^{1}\log(\sqrt{1-x^2})~p_n^2(x)~w(x)~dx+o(1)\,, \quad
n\to \infty\,.
\end{split}
\end{equation}
On the other hand, if for $\varepsilon \in(0,1)$ we denote
$$
\log_\varepsilon (x)\isdef \max \{ \log(\varepsilon  ),
\log(\sqrt{1-x^2}) \} \in C[-1,1]\,,
$$
then
\begin{align*}
& \left| \int_{-1}^{1}\log(\sqrt{1-x^2})~p_n^2(x)~w(x)~dx
-\frac{1}{\pi}\,
\int_{-1}^{1}\log(\sqrt{1-x^2})~\frac{dx}{\sqrt{1-x^2}}\right| \\
 \leq & \left|
\int_{-1}^{1}   \left( \log(\sqrt{1-x^2}) - \log_\varepsilon
(x)\right) \, p_n^2(x)~w(x)~dx \right| \\
+& \left| \int_{-1}^{1}\log_\varepsilon (x) \, \left( p_n^2(x)~w(x)-
\frac{1}{\pi \sqrt{1-x^2}}\right) \, dx \right| \\ +& \left|
\int_{-1}^{1} \left( \log_\varepsilon (x)- \log(\sqrt{1-x^2})
\right) \, \frac{1}{\pi \sqrt{1-x^2}} \, dx \right|=I_1+I_2+I_3\,.
\end{align*}
By \eqref{convErdos}, $I_2=o(1)$, as $n\to \infty$, while
$$
I_3=-\int_{\sqrt{1-x^2 }<\varepsilon } \log(\sqrt{1-x^2}) \,
\frac{1}{\pi \sqrt{1-x^2}} \, dx=\mathcal O(\varepsilon )\,.
$$
Using the same arguments as for \eqref{limitFirstIntegral} we find
also that $I_1=\mathcal O(\varepsilon )$. Taking into account that
$\varepsilon >0$ is arbitrary and using \eqref{limitFirstIntegral}
we obtain that
\begin{equation*}
    \begin{split}
\int_{-\alpha_n}^{\alpha_n}\log(\sqrt{\alpha_n^2-x^2})~p_n^2(x)~w(x)~dx
\\ =\frac{1}{\pi}\,
\int_{-1}^{1}\log(\sqrt{1-x^2})~\frac{dx}{\sqrt{1-x^2}}+o(1)\,,
\quad n\to \infty\,.
    \end{split}
\end{equation*}
The identity
$$
\frac{1}{\pi}\,
\int_{-1}^{1}\log(\sqrt{1-x^2})~\frac{dx}{\sqrt{1-x^2}}=-\log(2)
$$
is straightforward, which concludes the proof.
\end{proof}
\begin{remark}
This result is not surprising: if we denote by $\widehat \nu_n$ the
absolutely continuous measure on $[-\alpha _n, \alpha _n]$ with
$\widehat \nu'_n(x)=p_n^2(x)w(x)$, then $\widehat \nu_n \to \mu$ in
the weak-$*$ topology, where $\mu$ is the Robin measure of $\Delta$.
The integral in the left hand side of \eqref{limitPotentials} can be
rewritten as
$$
-\frac{1}{2}\,\left( V^{\widehat \nu_n}(-\alpha _n)+V^{\widehat
\nu_n}( \alpha _n)\right) \longrightarrow -\frac{1}{2}\,\left(
V^{\mu}(-1)+V^{\mu}( 1)\right)=-\log(2)\,.
$$
\end{remark}
\medskip
Now we turn to the proof of Theorem \ref{thFn}. Using function
$\Phi_n$ introduced in \eqref{defPhi}, let us denote
\begin{equation}
\begin{split}
 f_n(\alpha_n\cos\theta)& \isdef
\sqrt{\alpha_n}~p_n(\alpha_n\cos\theta)~\sqrt{w(\alpha_n\cos\theta)}~\sqrt{\sin\theta}
,\label{deffn}\\
g_n(\alpha_n\cos\theta)& \isdef
\sqrt{\frac{2}{\pi}}\cos\left(\frac{\theta}{2}-\frac{\pi}{4} +n
\Phi_n(\theta)\right).
\end{split}
\end{equation}
Since $w\in\mathcal F(C^2+)$, by \cite[Theorem 15.1 and Lemma
15.4]{Levin01},
\begin{equation}\label{asnorma1}
\int_0^\pi
 \big|f_n(\alpha_n\cos\theta)-g_n(\alpha_n\cos\theta)
\big|\, d\theta=o(1)\,, \quad n\to \infty\,,
 \end{equation}
(where we have used that on the bounded interval convergence in
$L^2$ is stronger than in $L^1$). Let us rewrite the definition of
$F_n$ as
\begin{align*}
F_n = & - \int_{-\alpha_n}^{\alpha_n} \log(p_n^2(x)~w(x)) p_n^2(x)
w(x)~dx\\
&
-\int_{\Delta\setminus[-\alpha_n,\alpha_n]}\log(p_n^2(x)~w(x))~p_n^2(x)~w(x)~dx
\\  =
&
-\int_{-\alpha_n}^{\alpha_n}\log(p_n^2(x)~w(x))~p_n^2(x)~w(x)~dx+o(1)\,,
\end{align*}
where we have used Lemma \ref{lemaFnMRS}. With the notation
\eqref{defR} and (\ref{deffn}) it is equivalent to
\begin{equation}\label{Fn2}
\begin{split}
F_n &= -\int_{-\alpha_n}^{\alpha_n} \mathcal{R}(f_n(x))~\frac{dx}{
\sqrt{\alpha_n^2-x^2}}+ \int_{-\alpha_n}^{\alpha_n}
\log(\sqrt{\alpha_n^2-x^2})~p_n^2(x)~w(x)~dx+o(1) \\ &=
-\int_{0}^{\pi}\mathcal{R}(f_n(\alpha_n\cos\theta))~d\theta-\log
2+o(1)\,, \quad n\to \infty\,,
\end{split}
\end{equation}
(see Lemma \ref{lemaFn}). Since by \cite[Theorem 1.17]{Levin01}
there exists a constant $M>\sqrt{2/\pi}$ such that for all $n\in
\N$, $ \big| f_n(\alpha_n\cos\theta)\big|\leq M$, for $
 \theta\in [0,\pi]$,  we get by (\ref{asnorma1})
\begin{align*}
&\int_0^\pi\left|\mathcal{R}(f_n(\alpha_n\cos\theta))-\mathcal{R}(g_n(\alpha_n\cos\theta))\right|d\theta\nonumber\\
\leq &
\max_{y\in[0,M]}\mathcal{R}'(y)\int_0^\pi\left|f_n(\alpha_n\cos\theta)-
g_n(\alpha_n\cos\theta)\right|d\theta=o(1)\,, 
\end{align*}
which yields
$$
F_n=-\int_{0}^{\pi}\mathcal{R}
\left(g_n(\alpha_n\cos\theta)\right)\,d\theta-\log 2+o(1)\,, \quad
n\to \infty\,.
$$

With the change of variable $\eta=\Phi_n(\theta)$ we rewrite
\begin{align}
 F_n &= -\int_{0}^{\pi}\mathcal{R}
\left(\sqrt{\frac{2}{\pi}}\cos\left(\Phi_n^{[-1]}( \eta)
-\frac{\pi}{4}+n\eta\right)\right)~ \frac{d\eta
}{\Phi_n'(\Phi_n^{[-1]}( \eta))} -\log
2+o(1) \nonumber\\
& =-\int_{0}^{\pi}\mathcal{R}
\left(\sqrt{\frac{2}{\pi}}\cos\left(\Phi_n^{[-1]}( \eta)
-\frac{\pi}{4}+n\eta\right)\right)~
d\eta-\log 2+o(1)\nonumber\\
& =-\int_{0}^{\pi}\mathcal{R}
\left(\sqrt{\frac{2}{\pi}}\cos\left(\eta-\frac{\pi}{4}+n\eta\right)\right)~
d\eta-\log 2+o(1)\label{lastStepForF},
\end{align}
where we have used Lemma \ref{lmPhin}. It remains to use the
following analogue of the Lebesgue lemma, proved in
\cite{Aptekarev:95} under weaker conditions:
\begin{lm}\label{lmAptekarev}
Let $g$ be a $\pi$-periodic continuous function on $[0,+\infty)$,
and $h\in C[0,\pi]$. Then
$$\int_0^\pi
f(n\theta+h(\theta))~d\theta=\int_0^\pi f(\theta)d\theta+o(1).$$
\end{lm}
Applying this lemma to \eqref{lastStepForF}, we obtain finally
$$
F_n=-1+\log(2)+\log(\pi)-\log(2)+o(1)=\log(\pi)-1+o(1)\,, \quad
n\to \infty\,.
$$
\section{Asymptotics of $G_n$: proof of Theorem \ref{thGn}}\label{SecGn}
We start again with some technical results:
\begin{lm}\label{lmGn2}
When $n\rightarrow\infty$,
$$\int_{0}^1
x(1-x^2)~s'(x)~p_n^2(x)~w(x)dx=B_2+o(1),$$ where
\begin{equation}\label{defB3}
B_2=\int_{0}^1x(1-x^2)~s'(x)\frac{1}{\pi \sqrt{1-x^2}}~dx\,.
\end{equation}
\end{lm}
\begin{proof} By \eqref{convErdos}, it is sufficient to show that
$(1-x^2) s'(x)$ can be extended as a continuous function to the
whole interval $\Delta$; for this purpose we only need to show that
the limit
$$\lim_{x\rightarrow 1^-}(1-x^2)~s'(x) $$
exists. From the explicit expression for $s$
 it is easy to find that for $x\in (0,1)$,
\begin{align*}
(1-x^2) s'(x) =& (1-x^2)\, \frac{w'(x)}{w(x)}-(1-x^2)\,
\frac{w_0'(x)}{w_0(x)} =(1-x^2)~\frac{w'(x)}{w(x)}+\frac{2\pi
a}{\sqrt{1-x^2}}\\ =& -2x(\lambda-1/2)-2ax+
 2 \arccos x~\frac{a}{\sqrt{1-x^2}} \\ &  -2\frac{t}{x}
  \left(\Im\psi(\lambda+i
  t)-\frac{\pi}{2}\right).
\end{align*}
It remains to use \eqref{limitImPsi}, and the statement follows.
\end{proof}
Let us denote
\begin{equation}\label{defbeta} p_n(x)=\gamma_n
x^n+\beta_n x^{n-2}+\text{lower degree terms};
\end{equation}
the explicit expression for $\gamma_n$ was given in
\eqref{leadingCoeff}.
\begin{lm}\label{lmGn3}
The followings identities hold:
\begin{equation}\label{eq1lm3}
\int_{-1}^1x~p_n(x)~p_n'(x)~w(x)~dx=n\,,
\end{equation}
and
\begin{equation}\label{eq2lm3}
\int_{-1}^1x^3 p_n(x)~p_n'(x)~w(x)~dx =
n(a_{n+1}^2+a_{n}^2)-2a_na_{n-1}\frac{\beta_n}{\gamma_{n-2}},
\end{equation}
where $a_n$ are the coefficients of the recurrence relation
(\ref{RecPoll}), and $\gamma_n$, $\beta_n$ are the coefficients of
$p_n$ defined in (\ref{defbeta}).
\end{lm}
\begin{proof}
By the recurrence relation (\ref{RecPoll}),
$$xp_n(x)p_n'(x)=a_{n+1}p_{n+1}(x)p_n'(x)+a_np_{n-1}(x)p_n'(x),$$
so that
\begin{align*}
\int_{-1}^1xp_n(x)p_n'(x)~w(x)~dx&=a_n\int_{-1}^1p_{n-1}(x)p_n'(x)~w(x)~dx\\
&=a_n \frac{n \gamma_n}{\gamma_{n-1}}
\int_{-1}^1p_{n-1}^2(x)~w(x)~dx=n,
\end{align*}
where we have used the well known fact that
$a_n=\gamma_{n-1}/\gamma_n$. This proves (\ref{eq1lm3}).
Again, from (\ref{RecPoll}) it is easy to find that
$$x^2p_n(x)=a_{n+2}a_{n+1}p_{n+2}(x)+(a_{n+1}^2+a_n^2)p_n(x)+a_{n}a_{n-1}p_{n-2},$$
and we get
\begin{align}\label{dmlm3-1}
\int_{-1}^1x^3&p_n(x)p_n'(x)~w(x)~dx
=(a_{n+1}^2+a_n^2)\int_{-1}^1xp_n(x)p_n'(x)~w(x)~dx \nonumber\\& +
a_na_{n-1}\int_{-1}^1xp_{n-2}p_n'(x)~w(x)~dx.
\end{align}
First integral in the right hand side was computed in
(\ref{eq1lm3}), and it remains to concentrate our attention on the
second one. Since
\begin{align*}
xp_n'(x)&=n \gamma_n x^n+(n-2)\beta_nx^{n-2}+\text{ lower degree
terms }\\
&=np_n(x)-2\frac{\beta_n}{\gamma_{n-2}}p_{n-2}(x)+\text{ lower
degree terms,}
\end{align*}
we have
$$\int_{-1}^1xp_{n-2}(x)p_n'(x)~w(x)~dx=-2\frac{\beta_n}{\gamma_{n-2}}.$$
Substituting it in (\ref{dmlm3-1}), we obtain (\ref{eq2lm3}).
\end{proof}
We find next an expression for the ratio $\beta_n/\gamma_n$ in terms
of the coefficients of the recurrence relation:
\begin{lm}
With the notations introduced in (\ref{RecPoll}) and
(\ref{defbeta}),
\begin{equation}\label{eq1lm1}
\frac{\beta_{n+1}}{\gamma_{n+1}}=\frac{\beta_n}{\gamma_n}-a_n^2,\qquad
n\in\mathbb{N};
\end{equation}
in particular
\begin{equation}\label{eq2lm1}
\frac{\beta_{n+1}}{\gamma_{n+1}}=-\sum_{k=1}^n a_k^2.
\end{equation}
\end{lm}
\begin{proof}
Comparing the coefficients of $x^{n-1}$ in both sides of
(\ref{RecPoll}) we obtain that
$$\beta_n=a_{n+1}\beta_{n+1}+a_n \gamma_{n-1},$$
so that
$$\frac{\beta_{n+1}}{\gamma_{n+1}}=\frac{1}{a_{n+1}}~\frac{\gamma_{n}}{\gamma_{n+1}}~
\frac{\beta_n}{\gamma_n}-\frac{a_n}{a_{n+1}}~\frac{\gamma_{n-1}}{\gamma_{n}}~\frac{\gamma_{n}}{\gamma_{n+1}},$$
and so, using the identity
$$a_n=\frac{\gamma_{n-1}}{\gamma_n},$$
(\ref{eq1lm1}) holds.
Formula (\ref{eq2lm1}) follows from (\ref{eq1lm1}) and the fact
that $\beta_1=0$.
\end{proof}
\begin{remark}
Observe that we have used only the symmetry of the recurrence
relation, so these formulas are valid for any even weight function
on $[-1,1]$.
\end{remark}
 \begin{cor}\label{cor2}
 For the symmetric Pollaczek polynomials the following asymptotic formula is valid:
 \begin{equation}\label{cocbeta}
\frac{\beta_n}{\gamma_n}=-~\frac{n}{4}+\frac{a}{2} \log
 n-\frac{a-\lambda}{4}-\frac{a}{2}\psi(a+\lambda)+\mathcal{O}\left(\frac{1}{n}\right),\qquad n\rightarrow \infty.
 \end{equation}
 \end{cor}
\begin{proof}
From (\ref{CoefRecPoll}) it is easy to obtain that
$$a_k^2=1+\frac{a^2-\lambda^2-a+\lambda}{k+\lambda+a-1}+\frac{\lambda^2-a^2-\lambda-a}{k+\lambda+a},$$
and (\ref{eq2lm1}) lead us to
\begin{align*}
\frac{\beta_n}{\gamma_n}= &
-\frac{n-1}{4}+\frac{1}{4}\sum_{k=1}^{n-1}\frac{2a}{k+\lambda+a-1}
+\frac{1}{4}\frac{-\lambda^2+a^2+\lambda+a}{n+\lambda+a-1}
\\ & -\frac{1}{4}\frac{-\lambda^2+a^2+\lambda+a}{\lambda+a}.
\end{align*}
Using that
\begin{equation*}
\sum_{k=1}^{n-1}\frac{1}{k+\lambda+a-1}=\psi(n+\lambda+a-1)-\psi(\lambda+a)=\log(n)-\psi(\lambda+a)+
\mathcal{O}\left(\frac{1}{n}\right),
\end{equation*}
when $n\to\infty$, we finally get
$$\frac{\beta_n}{\gamma_n}=-\frac{n-1}{4}+\frac{2a}{4}(\log(n)-\psi(\lambda+a))
+\frac{1}{4}(\lambda-a-1)+\mathcal{O}\left(\frac{1}{n}\right),$$
which is equivalent to the statement of the Lemma.
\end{proof}
Now we can prove Theorem \ref{thGn}. Remember that the technique of
\cite{Levin03} is not valid here because an additional assumption on
$w$ from \cite{Levin03} is not satisfied. The central idea in our
proof is to take advantage of the fact that the main contribution to
the asymptotics of $G_n$ comes from the behavior of the weight $w$
at the endpoints of $\Delta$ (see Section \ref{SecPeso}). Using
functions $w_0$ and $s$ introduced in (\ref{defw0}), we write $G_n$
in the form
\begin{align}
G_n&=\int_{-1}^1
\log(w_0(x))~p_n^2(x)~e^{s(x)}~w_0(x)~dx+\int_{-1}^1
s(x)~p_n^2(x)~w(x)~dx.\nonumber
\end{align}
In particular, since $s\in C[-1,1]$ (see Lemma \ref{prolims}),
applying (\ref{convErdos}) in the second integral we have
\begin{equation}\label{Gn}
G_n=\int_{-1}^1 \log(w_0(x))~p_n^2(x)~e^{s(x)}~w_0(x)~dx+B_1+o(1),
\end{equation}
where
\begin{align*}
B_1&=\lim_{n\rightarrow\infty}\int_{-1}^1
s(x)~p_n^2(x)~w(x)~dx=\int_{-1}^1
\frac{s(x)}{\pi~\sqrt{1-x^2}}~dx.
\end{align*}
If we denote
$$g(x)=\log(w_0(x))=-~\frac{2\pi a |x|}{\sqrt{1-x^2}},$$
taking into account the symmetry, we can rewrite the integral in
the right hand side of (\ref{Gn}) as
\begin{align*}
&\int_{-1}^1
g(x)~p_n^2(x)~e^{s(x)}~e^{g(x)}~dx=2\int_0^1\frac{g(x)}{g'(x)}~p_n^2(x)~e^{s(x)}~e^{g(x)}g'(x)~dx
\\&\quad =2\int_0^1\left(\frac{g(x)}{g'(x)}~p_n^2(x)~e^{s(x)}\right)~de^{g(x)}.
\end{align*}
Observe that  for $x \in [0,1]$, $g(x)/g'(x)=x(1-x^2)$, so
integrating by parts,
\begin{align*}
&\int_{-1}^1
g(x)~p_n^2(x)~e^{s(x)}~e^{g(x)}~dx\\&\qquad=2\left[x(1-x^2)~p_n^2(x)~w(x)\right]_{x=0}^{x=1}-
2\int_0^1\left(x(1-x^2)~p_n^2(x)~e^{s(x)}\right)'~w_0(x)~dx\\
&\qquad=-\int_{-1}^1(1-3x^2)p_n^2(x)~w(x)~dx-\int_{-1}^1x(1-x^2)~\left(p_n^2(x)\right)'~w(x)~dx\\&\qquad\quad-
2\int_0^1x(1-x^2)~p_n^2(x)~s'(x)~w(x)~dx\,.
\end{align*}
The asymptotics of each of these three integrals can be computed by
means of (\ref{convErdos}), Lemma \ref{lmGn3}, and Lemma
\ref{lmGn2}, respectively, obtaining that
\begin{align*}
G_n&=2n(a_{n+1}^2+a_n^2-1)-4a_na_{n-1}\frac{\beta_n}{\gamma_{n-2}}+B_1-2B_2+\frac{1}{2}+o(1)\\
&=2n(a_{n+1}^2+a_n^2-1)-4\frac{\beta_n}{\gamma_n}+B_1-2B_2+\frac{1}{2}+o(1),
\quad n \to \infty\,.
\end{align*}
Using that $a_n\rightarrow 1/2$ (see \eqref{CoefRecPoll}) and
\eqref{cocbeta}, we get that
\begin{equation}\label{GnPreliminary}
G_n=-2a \log(n)+B_1-2B_2+\frac{1}{2}+B_3+o(1),
\end{equation}
where
\begin{equation}\label{B3}
B_3=-2a+(a-\lambda)+2a\psi(a+\lambda)=-a-\lambda+2a\,
\psi(a+\lambda).
\end{equation}
Let us simplify the expression of the constant term of this
asymptotics. First,
\begin{equation*}
B_1-2B_2=\int_{-1}^1
\frac{s(x)}{\pi~\sqrt{1-x^2}}~dx-2\int_{0}^1x\sqrt{1-x^2}~s'(x)\frac{1}{\pi
}~dx,
\end{equation*}
and integrating by parts the second integral,
\begin{align}
B_1-2B_2&= 2\int_{0}^1
\frac{s(x)}{\pi~\sqrt{1-x^2}}-\frac{2}{\pi}\left[x\sqrt{1-x^2}s(x)\right]_{x=0}^{x=1}
+2\int_0^1\frac{1-2x^2}{\pi\sqrt{1-x^2}}~s(x)~dx\nonumber\\
&=\frac{4}{\pi}\int_0^1s(x)~\sqrt{1-x^2}~dx.\label{B12}
\end{align}
Using  the explicit expression for $s$ on $[0,1]$, the right hand
side in (\ref{B12}) is reduced to
\begin{align}
&\quad\log\left(\frac{2^{2\lambda}~(\lambda+a)}{2
\pi~\Gamma(2\lambda)}\right)~\frac{4}{\pi}\int_0^1\sqrt{1-x^2}dx+
\frac{4(\lambda-1/2)}{\pi}\int_0^1\log(1-x^2)~\sqrt{1-x^2}~dx \nonumber\\
&\quad +\frac{4}{\pi}\int_0^1 2ax\arccos
x~dx+\frac{4}{\pi}\int_0^1\log\left|\Gamma(\lambda+i
 t)\right|^2~\sqrt{1-x^2}~dx+4\int_0^1ax~dx  \nonumber \\
 &=\log\left(\frac{2^{2\lambda}~(\lambda+a)}{2
\pi~\Gamma(2\lambda)}\right)+(1-2\log(2))(\lambda-1/2)  \nonumber  \\
&\quad +a+\frac{4}{\pi}\int_0^1\log\left|\Gamma(\lambda+i
 t)\right|^2~\sqrt{1-x^2}~dx+2a  \nonumber \\
 &=\log\left(\frac{(\lambda+a)}{
\pi~\Gamma(2\lambda)}\right)+\lambda+3a+\frac{4}{\pi}\int_0^1\log\left|\Gamma\left(\lambda+i
 \frac{ax}{\sqrt{1-x^2}}\right)\right|^2~\sqrt{1-x^2}~dx.
 \label{B2B3}
 \end{align}
Let us compute now the value of this last integral. With the change
of variables $u=x/\sqrt{1-x^2}$ we obtain that
\begin{align*}
\int_0^1\log\left|\Gamma\left(\lambda+i
 \frac{ax}{\sqrt{1-x^2}}\right)\right|^2~\sqrt{1-x^2}~dx & =\int_0^{+\infty}\log\left|\Gamma\left(\lambda+i
a u\right)\right|^2~\frac{du}{(1+u^2)^2}\\
& =\int_{-\infty}^{+\infty}\log\left|\Gamma\left(\lambda+i a
u\right)\right|~\frac{du}{(1+u^2)^2}\,.
\end{align*}
With $\lambda >0$, $a\geq 0$, function
$$
f(u)\isdef  \frac{\log\left(\Gamma\left(\lambda+i a
u\right)\right)}{(1+u^2)^2}
$$
is meromorphic and single valued in the lower half plane $\{ \Im (
u) <0\}$, with a double pole at $u=-i$. Taking into account that by
Stirling formula,
$$
\log\left(\Gamma\left(\lambda+i a u\right)\right) \sim  (\lambda+i a
u) |u| \log(\lambda+i a u) \quad \text{as } u \to \infty, \; \Im (
u) <0\,,$$
we may apply the residue calculus to establish that
$$
\int_{-\infty}^{+\infty}\log\left( \Gamma\left(\lambda+i a
u\right)\right)~\frac{du}{(1+u^2)^2}=-2\pi i \res_{u=-i}
f(u)=\frac{\pi}{2}\, \left(a\psi (\lambda +a)-\log (\Gamma(\lambda
+a)) \right)\,.
$$
Taking the real part, we get that
\begin{equation}\label{valueIntegral}
 \frac{4}{\pi}
\int_0^1\log\left|\Gamma\left(\lambda+i
 \frac{ax}{\sqrt{1-x^2}}\right)\right|^2~\sqrt{1-x^2}~dx=
 2 a\psi (\lambda +a)- 2 \log (\Gamma(\lambda
+a))  \,.
\end{equation}
Gathering \eqref{B3}--\eqref{valueIntegral} in
\eqref{GnPreliminary} we conclude the proof of Theorem \ref{thGn}.

\medskip

\begin{remark} The idea of this proof can be applied also to the case of a
non-symmetric weight of the form
\begin{equation*}\label{descomp3}
w(x)=\exp\left\{-\frac{4c}{(1-x)^{\alpha}}-\frac{4d}{(1+x)^{\alpha}}+s(x)\right\},
\end{equation*}
where $s\in C^1 [-1,1] $, and  $\alpha \in [1/2,1]$. 
\end{remark}

\section{Proof of corollaries \ref{CorAsPoll} and \ref{CorEnergy}}\label{SecCorLevin}

\emph{Proof of Corollary \ref{CorAsPoll}:} Consider the integral
$$I_n=-\frac{2}{\pi}\int_{\alpha_{-n}}^{\alpha_n}\frac{Q(x)}
{\sqrt{(\alpha_n-x)(x-\alpha_{-n})}}~dx=-\frac{4}{\pi}\int_{0}^{\alpha_n}\frac{Q(x)}
{\sqrt{(\alpha_n-x)(x-\alpha_{-n})}}~dx.
$$
By (\ref{defw0}),
\begin{align*}
I_n = & -\frac{4}{\pi}\int_{0}^{\alpha_n}\left(\frac{\pi
a|x|}{\sqrt{1-x^2}}+\frac{1}{2}\log(w(0))-\frac{s(x)}{2}\right)
~\frac{1}{\sqrt{\alpha_n^2-x^2}}~dx\\
 = & -\frac{4}{\pi}\int_{0}^{1}\left(\frac{\pi a|\alpha_n
x|}{\sqrt{1-\alpha_n^2x^2}}
+\frac{1}{2}\log(w(0))-\frac{s(\alpha_nx)}{2}\right)
~\frac{1}{\sqrt{1-x^2}}~dx\nonumber\\
 = & -\frac{4}{\pi}\int_{0}^{1}\frac{\pi a|\alpha_n
x|}{\sqrt{1-\alpha_n^2x^2}}~\frac{1}{\sqrt{1-x^2}}~dx\\
& -\frac{4}{\pi}\int_{0}^{1}
\left(\frac{1}{2}\log(w(0))-\frac{s(\alpha_n
x)}{2}\right)~\frac{1}{\sqrt{1-x^2}}~dx.
\end{align*}
Since $s\in C[-1,1]$, the second integral is bounded; hence, by
\eqref{ellInt1},
$$I_n=-4 a \alpha_n~\frac{-1}{2}~\log(1-\alpha_n)+\mathcal O(1)=2a\log(1-\alpha_n)+\mathcal O(1)\,,
\quad n \to \infty\,.$$ Finally, from the asymptotics of $\alpha_n$
found in Proposition \ref{proMRS} we obtain
$$I_n=2a\log(a/n)+\mathcal O(1)=-2a\log (n)+\mathcal O(1)\,,
\quad n \to \infty\,,$$ and comparing this expression with the
result of Corollary \ref{TeoAsPoll}, the statement follows.

\medskip

\emph{Proof of Corollary \ref{CorEnergy}:} Taking into account the
relation between the entropy and the mutual energy (\ref{kullback2})
we obtain that
\begin{equation*}
I[\rho _n, \nu_n]=\frac{E_n+2\log(\gamma_n)}{2n}\,.
\end{equation*}
By \eqref{leadingCoeff},
$$2 \log (\gamma_n)=2n\log(2)+2a\log(n)-\log\left(\frac{\Gamma(\lambda+a+1)
~\Gamma(\lambda+a)}{\Gamma(2\lambda)}\right)+o(1),
$$
so that by \eqref{asymptoticFormulaforE},
\begin{equation*}
I[\nu_n,\lambda_n]=\log(2)+\frac{1}{2n}\left(\tau(\lambda,a)-
\log\left(\frac{\Gamma(\lambda+a+1)~\Gamma(\lambda+a)}{\Gamma(2\lambda)}\right)\right)+o\left(\frac{1}{n}\right).
\end{equation*}
The use of the explicit expression for $\tau(\lambda,a)$ in
\eqref{def_tau} concludes the proof of the Corollary.

\section*{Acknowledgement}

The authors were supported, in part, by a research grant from the
Ministry of Science and Technology (MCYT) of Spain, project code
BFM2001-3878-C02, and by Junta de Andaluc\'{\i}a, Grupo de Investigaci\'{o}n
FQM229. The research of A.M.F. was supported also by Research
Network ``Network on Constructive Complex Approximation (NeCCA)'',
INTAS 03-51-6637, and by NATO Collaborative Linkage Grant
``Orthogonal Polynomials: Theory, Applications and
Generalizations'', ref. PST.CLG.979738. We are indebted also to
Prof. D.\ S.\ Lubinsky for very interesting discussions.


\begin{thebibliography}{10}
\bibitem{abramowitz/stegun:1972}
{\sc M.~Abramowitz and I.~A. Stegun}, {\em Handbook of Mathematical
Functions},
  Dover Publ., New York, 1972.
\bibitem{Aptekarev:95}
{\sc A.~I. Aptekarev, V.~S. Buyarov, and J.~S. Dehesa}, {\em
Asymptotic
  behavior of the {$L^p$}-norms and the entropy for general orthogonal
  polynomials}, Russian Acad.\ Sci.\ Sb.\ Math., 82 (1995), pp.~373--395.
\bibitem{Beckermann04}
{\sc B.~Beckermann, A.~Mart\'{\i}nez-Finkelshtein, E.~A.
Rakhmanov, and
  F.~Wielonsky}, {\em Asymptotic upper bounds for the entropy of orthogonal
  polynomials in the {Szeg\H{o}} class}, J. Math. Physics, 45 (2004),
  pp.~4239--4254.
\bibitem{Bialynicki-Birula:75}
{\sc I.~Bialynicki-Birula and J.~Mycielsky}, {\em Uncertainty
relations for
  information entropy in wave mechanics}, Commun.\ Math.\ Phys., 44 (1975),
  pp.~129--132.
\bibitem{Buyarov:97}
{\sc V.~S. Buyarov}, {\em On information entropy of {G}engenbauer
polynomials},
  Vesntik Moskow Univ., Ser. \textbf{1} 6 (1997), pp. 8--11
\newblock (in russian).
\bibitem{Buyarov:99a}
{\sc V.~S. Buyarov, J.~S. Dehesa, A.~Mart\'{\i}nez-Finkelshtein,
and E.~B.
  Saff}, {\em Asymptotics of the information entropy for {J}acobi and
  {L}aguerre polynomials with varying weights}, J. Approx. Theory, 99 (1999),
  pp.~153--166.
\bibitem{MR1790053}
{\sc V.~S. Buyarov, P.~L\'{o}pez-Art\'{e}s,
A.~Mart\'{\i}nez-Finkelshtein,
  and W.~Van~Assche}, {\em Information entropy of {G}egenbauer polynomials}, J.
  Phys. A, 33 (2000), pp.~6549--6560.
\bibitem{Buyarov04}
{\sc V.~Buyarov, J.~S. Dehesa, A.~Mart\'{\i}nez-Finkelshtein, and
  J.~S\'{a}nchez-Lara}, {\em Computation of the entropy of polynomials
  orthogonal on an interval}, SIAM J. Sci. Comp., 26 (2004), pp.~488--509.
\bibitem{Dehesa:97}
{\sc J.~S. Dehesa, W.~Van~Assche, R.~J. Y\'{a}\~{n}ez}, {\em
Information
  entropy of classical orthogonal polynomials and their application to the
  harmonic oscillator and {C}oulomb potentials}, Methods Appl.\ Anal. \textbf{4}
  (1997), pp. 91--110.
\bibitem{Dehesa:01}
{\sc J.~S. Dehesa, A.~Mart\'{\i}nez-Finkelshtein, and
J.~S\'{a}nchez-Ruiz},
  {\em Quantum information entropies and orthogonal polynomials}, J.\ Comput.\
  Appl.\ Math., 133 (2001), pp.~23--46.
\bibitem{Dehesa:98}
{\sc J.~S. Dehesa, R.~J. Ya\~{n}ez, A.~I. Aptekarev, and V.~S.
Buyarov}, {\em
  Strong asymptotics of {L}aguerre polynomials and information entropies of
  {2D} harmonic oscillator and {1D} {C}oulomb potentials}, J.\ Math.\ Physics,
  39 (1998), pp.~3050--3060.
\bibitem{Gradshtein95}
{\sc I.~S. Gradshtein and I.~M. Ryzhik}, {\em Table of Integrals,
Series and
  Products}, Academic Press, San Diego, CA, fifth~ed., 1995.
\bibitem{Levin01}
{\sc E. Levin and D.~S. Lubinsky}, {\em Orthogonal Polynomials for
  Exponential Weights}, vol.~4 of CMS Books in Mathematics, Springer Verlag,
  2001.
\bibitem{Levin03}
{\sc E.~Levin and D.~S. Lubinsky}, {\em Asymptotics for entropy
integrals
  associated with exponential weights}, J.\ Comput.\ Appl.\ Math., 156 (2003),
  pp.~265--283.
\bibitem{Luke75}
{\sc Y.~L. Luke}, {\em Mathematical Functions and their
Approximations},
  Academic Press, New York, 1975.
\bibitem{Mate}
{\sc A. M\'{a}t\'{e}, P. Nevai, V. Totik,} {\em Strong and weak
convergence of orthogonal polynomials}, Am. J. Math. {\bf 109}
(1987), pp. 239--282.
\bibitem{Parr89}
{\sc R.~G. Parr and W.~Yang}, {\em Density Functional Theory of
Atoms and
  Molecules}, Oxford University Press, New York, 1989.
\bibitem{Rakhmanov:77}
{\sc E.~A. Rakhmanov}, {\em On the asymptotics of the ratio of
orthogonal
  polynomials}, Math. USSR Sb., 32 (1977), pp.~199--213.
\bibitem{Saff:97}
{\sc E.~B. Saff and V.~Totik}, {\em Logarithmic Potentials with
External
  Fields}, vol.~316 of Grundlehren der Mathematischen Wissenschaften,
  Springer-Verlag, Berlin, 1997.
\bibitem{Lara02}
{\sc J.~F. S\'{a}nchez~Lara}, {\em On the asymptotic expansion of
the entropy
  of {G}egenbauer polynomials}, J.\ Comput.\ Appl.\ Math., 142 (2002),
  pp.~401--409.
\bibitem{szego:1975}
{\sc G.~Szeg\H{o}}, {\em Orthogonal Polynomials}, vol.~23 of Amer.\
Math.\
  Soc.\ Colloq.\ Publ., Amer.\ Math.\ Soc., Providence, {RI}, fourth~ed., 1975.
\bibitem{Yanez:94}
{\sc R.~J. Y\'{a}\~{n}ez, W.~Van~Assche, and J.~S. Dehesa}, {\em
Position and
  momentum information entropies of the {D}-dimensional harmonic oscillator and
  hydrogen atom}, Physical Rev. {A}, 50 (1994), pp.~3065--3079.
  \bibitem{Yanez:99}
{\sc R.~J. Ya\~{n}ez, W.~Van~Assche, R.~Gonz\'{a}lez-F\'{e}rez, J.~S.
  Dehesa}, {\em Entropic integrals of hyperspherical harmonics and spatial
  entropy of {$D$}-dimensional central potentials}, J. Math. Phys. \textbf{40}
  (1999), pp. 5675--5686.
\end{thebibliography}

\end{document}